\documentclass[11pt]{article}
\usepackage{latexsym,amsmath,color,amsthm,amssymb,epsfig,graphicx,mathrsfs}
\usepackage{graphicx}
\usepackage{amssymb}
\usepackage[left=0.9in,top=0.9in,right=0.9in,bottom=0.9in]{geometry}
\usepackage[linktocpage=true]{hyperref}
\usepackage{setspace}
\usepackage{amssymb, amsmath, amsthm, graphicx,mathrsfs}
\usepackage{caption}

\usepackage{comment}
\usepackage{tikz}
\def\qed{\hfill\ifhmode\unskip\nobreak\fi\quad\ifmmode\Box\else\hfill$\Box$\fi}
\def\ite#1{\hfill\break${}$\hbox to 50pt {\quad(#1)\hfill}}
\newtheorem{thm}{Theorem}[section]

\newtheorem{lem}[thm]{Lemma}
\newtheorem{conj}[thm]{Conjecture}

\newtheorem{prop}{Proposition}[section]

\def\ex{{\rm{ex}}}
\def\mc{{d_{\min}}}

\begin{document}

\title{ Hypergraphs not containing a tight tree with a bounded trunk ~II:
3-trees with a trunk of size 2
}

\pagestyle{myheadings} \markright{{\small{\sc F\"uredi, Jiang, Kostochka, Mubayi, Verstra\"ete:  Tur\'an numbers of tight trees 
}}}

\author{
\hspace{0.8in} Zolt\'an F\" uredi\thanks{Research supported by grant K116769
from the National Research, Development and Innovation Office NKFIH, and
by the Simons Foundation Collaboration grant \#317487.}
\and
Tao Jiang\thanks{Research partially supported by NSF award DMS-1400249.}
\and
Alexandr Kostochka\thanks{Research of this author is supported in part by NSF grant
 DMS-1600592 and by grants  18-01-00353A and 16-01-00499
  of the Russian Foundation for Basic Research.
} \hspace{0.8in} \and
Dhruv Mubayi\thanks{Research partially supported by NSF award DMS-1300138.} \and Jacques Verstra\"ete\thanks{Research supported by NSF award DMS-1556524.}
}

\date{ July 16,  2018} 

\maketitle

\vspace{-0.3in}

\begin{abstract}
A {\it tight $r$-tree} $T$ is an $r$-uniform hypergraph that has an edge-ordering
$e_1, e_2, \dots, e_t$ such that for each $i\geq 2$, 
 $e_i$ has a vertex $v_i$ that does not belong to any previous edge and  $e_i-v_i$ is contained in $e_j$ for some $j<i$.
 Kalai~conjectured in 1984  that every  $n$-vertex  $r$-uniform hypergraph with more than
 $\frac{t-1}{r}\binom{n}{r-1}$ edges contains every tight $r$-tree $T$ with $t$ edges.

A {\em trunk} $T'$ of a tight $r$-tree $T$ is a tight subtree $T'$ of $T$ such that vertices in $V(T)\setminus V(T')$ 
are leaves in $T$.
Kalai's Conjecture was proved in 1987 for tight $r$-trees that have a trunk of size one. 
In a previous paper we proved  an asymptotic version of Kalai's Conjecture for all tight $r$-trees that have
a trunk of bounded size.  In this paper we continue that work 
to establish the exact form of Kalai's Conjecture for all tight $3$-trees with at least $20$ edges that have a trunk of size two.

\end{abstract}

\section{Introduction. Trees, trunks, and Kalai's conjecture
}
  For an $r$-uniform hypergraph ({\em $r$-graph}, for short) $H$, the {\em Tur\'an number} $\ex_r(n,H)$  is the largest
$m$ such that there exists an $n$-vertex $r$-graph $G$ with $m$ edges that does not contain $H$. 
Estimating $\ex_r(n,H)$ is a difficult problem
 even for $r$-graphs with a simple structure.
Here we consider Tur\'{a}n-type problems for so called tight $r$-trees. 
A {\it  tight $r$-tree} ($r\geq 2$)
is an $r$-graph 
whose edges can be ordered  so that 
each edge $e$ apart from the first one contains a vertex
$v_e$ that does not belong to any preceding edge
but the set $e-v_e$ is contained in some preceding edge. 
Such an ordering is called a {\em proper ordering} of the edges. 
A  usual graph tree is a tight $2$-tree. 

A vertex $v$ in a tight $r$-tree $T$ is a {\it leaf} if it has degree one in $T$.
A {\em trunk} $T'$ of a tight $r$-tree $T$ is a tight subtree of $T$ such that in some proper ordering of the edges of $T$ 
the edges of $T'$ are listed first and vertices in $V(T)\setminus V(T')$ are leaves in $T$.
Hence, each $e\in E(T)\setminus E(T')$ contains an $(r-1)$-subset of some $e'\in E(T')$
and a leaf in $T$ (that lies outside $V(T')$). In the case of $r=2$ each  $e\in E(T)\setminus E(T')$  is a pendant edge. 
Every tight tree $T$ with at least two edges has a trunk (for example, $T$ minus the last edge
in a proper ordering is a trunk).
Let $c(T)$ denote the minimum size of a trunk of $T$. 
We write $e(H)$ for the number of edges in $H$.

In this paper we consider the following classical conjecture.
\begin{conj}[Kalai 1984, see in~\cite{FF}] \label{kalai}
Let $T$ be a tight $r$-tree with $t$ edges.
Then $\ex_r(n,T)\leq \frac{t-1}{r}\binom{n}{r-1}$.
\end{conj}

The coefficient $(t-1)/r$ in this conjecture, if it is true,  is optimal as 
one can see using constructions obtained from partial Steiner systems due to R\"odl~\cite{Rodl}. 
The conjecture turns out to be  difficult even for very special cases of tight trees, in fact for $r=2$ it is the 
famous Erd\H os-S\'os conjecture. 
The following partial result on Kalai's conjecture was proved in 1987. 

\begin{thm} [
\cite{FF}] \label{FF}
Let $T$ be a tight $r$-tree with $t$ edges and $c(T)=1$.
Suppose that $G$ is an $n$-vertex $r$-graph with $e(G)>\frac{t-1}{r}\binom{n}{r-1}$. 
Then $G$ contains a copy of $T$. 
\end{thm}

In a previous paper~\cite{FJKMV2}, we showed that Conjecture \ref{kalai} holds {\em asymptotically} for tight $r$-trees with a
trunk of a bounded size. Our result is as follows. 
Define $a(r,c):=(r^r+1-\frac{1}{r})(c-1)$.
\begin{thm} [\cite{FJKMV2}] \label{main}
Let $T$ be a tight $r$-tree with $t$ edges and $c(T)\leq c$. 
Then
\[ 
\ex_r(n,T)\leq  \left(\frac{t-1}{r}+ a(r,c)\right)\binom{n}{r-1}.
   \]  
\end{thm}
The goal of this paper is to prove the conjecture in {\em exact} form
   for infinitely many $3$-trees.
\begin{thm} \label{two-central}
Let $T$ be a tight $3$-tree with $t$ edges and $c(T)\leq 2$. If $t\geq 20$ then  
 $$
    \ex_3(n,T)\leq \frac{t-1}{3} \binom{n}{2}.$$
\end{thm}

Beside ideas and observations from~\cite{FJKMV2}, discharging is quite helpful here.


\section{Notation and preliminaries. Shadows and default weights}
In this section, we introduce  some notation and list a couple of simple observations from~\cite{FJKMV2}.
For the sake of self-containment, we present their simple proofs as well. 

The {\it shadow} of an $r$-graph $G$ is
$\enskip\partial(G):=\{S: |S|=r-1,\;\mbox{\em and} \; S\subseteq e \;\mbox{\em for some}\;  e\in e(G)\}.$

 The {\it link} 
of  a set $D\subseteq V(G)$ in an $r$-graph $G$ is defined as 
$\enskip L_G(D):= \{ e\setminus D: e\in E(G), D\subseteq e\}.$

The {\em degree} of $D$, $d_G(D)$, is the number the edges of $G$  containing $D$.
If $G$ is an $r$-graph and $|D|=r-1$, the elements of $L_G(D)$ are vertices. In this case,
we also use $N_G(D)$ to denote $L_G(D)$. 
Many times we  drop the subscript $G$. 
For  $1\leq p\leq r-1$,  the {\it minimum $p$-degree} of $G$ is 
\[\delta_p(G):=\min\{d_G(D): |D|=p, \enskip \mbox{\em and}\enskip D\subseteq e \enskip \mbox{\em for some}\enskip  e\in E(G) \}.\]
For an $r$-graph $G$ and $D\in \partial(G)$, let $w(D):=\frac{1}{d_G(D)}$.
For each $e\in E(G)$, let 
\[ 
w(e):=\sum_{D\in \binom{e}{r-1}}  w(D)= \sum_{D\in \binom{e}{r-1}}\frac{1}{d_G(D)}.
\] 
We call $w$ the {\it default weight function} on $E(G)$ and $\partial(G)$. 
 Frankl and F\"uredi~\cite{FF} (and later some others) used the following simple property of this function.
 
\begin{prop} \label{weight}
Let $G$ be an $r$-graph. Let $w$ be the default weight function on $E(G)$ and $\partial(G)$. Then \[\sum_{e\in E(G)} w(e)=|\partial(G)|.\]
\end{prop}
\noindent
{\em Proof.}\enskip 
By definition, 
\[\sum_{e\in E(G)} w(e) =\sum_{e\in E(G)} \left( \sum_{D\in \binom{e}{r-1}} \frac{1}{d_G(D)}\right) 
=\sum_{D\in \partial(G)} \left(  \sum_{e\in E(G), D\subseteq e} \frac{1}{d_G(D)} \right)
=\sum_{D\in \partial(G)} 1 =|\partial(G)|.\qed\]

An {\it embedding} of an $r$-graph $H$ into an $r$-graph $G$ is an injection $f: V(H)\to V(G)$ such that
for each $e\in E(H)$, $f(e)\in E(G)$.
 The following proposition is folklore. 

\begin{prop} 
Let 
$G$ be an $r$-graph with $e(G)>q|\partial(G)|$. Then $G$ contains a subgraph $G'$ with $\delta_{r-1}(G')\geq \lfloor q\rfloor +1$.
\end{prop}
\noindent
{\em Proof.}\enskip 
Starting from $G$, if there exists $D\in \partial(G)$ of degree at most  $\lfloor q\rfloor$
in the current $r$-graph, we remove  the edges of this $r$-graph  containing $D$. Let $G'$ be
the final $r$-graph. Since we have deleted at most $q|\partial(G)|<e(G)$ edges,
$G'$ is nonempty. By the stopping rule, $\delta_{r-1}(G')\geq \lfloor q\rfloor +1$. \qed



\section{Lemmas for Theorem \ref{two-central}}

The idea behind the proof of  Theorem \ref{two-central} is
to find in the host $3$-graph $G$ a special pair of edges with good properties  where we plan to map the trunk of size~2 of $T$.
We use the weight argument together with discharging to find such special pairs in the next two lemmas.

Given edges $e=abc$ and $f=adc$ in a $3$-graph $G$ sharing pair $ac$, for a pair $\{x,y\}\subset \{a,b,c,d\}$, let
$d'_{e,f}(x,y)$ denote the number of $z\in V(G)\setminus \{a,b,c,d\}$ such that $xyz\in G$. By definition
\begin{equation}\label{d'}
d'_{e,f}(x,y)\geq d(x,y)-2\qquad\qquad \mbox{\em for every $\{x,y\}\subset \{a,b,c,d\}$. }
\end{equation}

\begin{lem} \label{discharging1} Let $m\geq 20$ be a positive integer and
let $G$ be a $3$-graph satisfying $e(G)>\frac{m}{3} |\partial(G)|$ and $\delta_2(G)>\frac{m}{3}$. Let $w$ be
the default weight function on $E(G)$ and $\partial(G)$.
Then there exist edges $e=abc$  and $f=adc$  in $G$ satisfying\\
{\rm (a)}\enskip  $w(e)<\frac{3}{m}$ and $w(ac)<\frac{1}{m}$,\\
{\rm (b)}\enskip  $\min\{d'_{e,f}(a,b),d'_{e,f}(c,b)\}\geq \left\lfloor\frac{m}{3}\right\rfloor$,\\
{\rm (c)}\enskip  $\max\{d'_{e,f}(a,b),d'_{e,f}(c,b)\}\geq \left\lfloor\frac{2m}{3}\right\rfloor$, and \\
{\rm (d)}\enskip  either $3(w(f)-\frac{3}{m})<(\frac{3}{m}-w(e))$ or   $\max\{d'_{e,f}(a,d),d'_{e,f}(c,d)\}\geq m-1$.
\end{lem}
\noindent
{\em Proof.}\enskip 
For convenience, let $w_0=\frac{3}{m}$. By Proposition \ref{weight}, $\sum_{e\in G} w(e)=|\partial(G)|$. So,
\begin{equation} \label{avg-wt}
\frac{1}{e(G)}\sum_{e\in G} w(e)=\frac{|\partial(G)|}{e(G)}<\frac{1}{m/3}=w_0.
\end{equation}
Hence the average weight of an edge in $G$ is less than $w_0$.
 We call an edge $e\in E(G)$ {\it light} if $w(e)<w_0$
and {\it heavy} otherwise. A pair $\{x,y\}$ of vertices in $G$ is {\em good}, if $d(xy)\geq m+1$.

To find the desired pair of edges $e,f$ we first do some marking of edges.
For every light edge $e$, fix an ordering, say $a,b,c$, of its vertices so
that $d(ab)\leq d(bc)\leq d(ac)$.  We call $ab, bc, ac$ the {\it low, medium, high} sides of $e$, respectively.

Since $e$ is light, $w(e)=\frac{1}{d(ab)}+\frac{1}{d(bc)}+\frac{1}{d(ac)}<w_0=\frac{3}{m}$, it follows that
\begin{equation} \label{weight-bounds}
d(ac)>m, \quad d(bc)>\frac{3m}{2}, \quad d(ab)>\frac{m}{3}.
\end{equation}
In particular, $ac$ is good. We define markings involving $e$ based on three cases.
 
{\em Case {\rm M1}: }  $d(ab)\geq \lfloor m/3\rfloor+2$ and $d(bc)\geq \lfloor 2m/3\rfloor+2$. In this case,
 we let $e$ {\em mark} every edge containing $ac$ apart from itself.

{\em Case {\rm M2}: } 
 $d(ab)\leq \lfloor m/3\rfloor+1$.  By \eqref{weight-bounds}, $d(ab)= \lfloor m/3\rfloor+1$, and since $e$ is light,
 \begin{equation}\label{m/3}
 d(ac)\geq d(bc)> \frac{1}{\frac{3}{m}-\frac{3}{m+3}}=\frac{m(m+3)}{9}.
 \end{equation} 
We let $e$ mark all the edges $acx\neq e$ containing $ac$ such that $abx$ is not an edge in $G$.
 By~(\ref{m/3}), in this case
  \begin{equation}\label{m/3'}
\mbox{\em $e$ marks at least} \quad \frac{m(m+3)}{9}-\frac{m+3}{3}=\frac{(m+3)(m-3)}{9}\quad\mbox{\em edges.}
 \end{equation}

{\em Case {\rm M3}: } 
$d(bc)\leq \lfloor 2m/3\rfloor+1$.  By \eqref{weight-bounds},  $d(bc)= \lfloor 2m/3\rfloor+1$.
Let $e$ mark all the edges $acx\neq e$ containing $ac$ such that $bcx$ is not an edge in $G$. 
Since $e$ is light, 
 \begin{equation}\label{2m/3}
 d(ac)> \frac{1}{\frac{3}{m}-2\frac{3}{2m+3}}=\frac{m(2m+3)}{9}.
 \end{equation} 
 Similarly to~(\ref{m/3'}), in this case
  \begin{equation}\label{2m/3'}
\mbox{\em $e$ marks at least} \quad \frac{m(2m+3)}{9}-\frac{2m+3}{3}=\frac{(2m+3)(m-3)}{9}\quad\mbox{\em edges.}
 \end{equation}  

We perform the above marking procedure for each light edge $e$.
\medskip

{\bf Claim 1.}  If $e$ is a light edge and $f$ is an edge marked by $e$ then (a)-(c) hold.
Further, if $f$ is light, then the lemma holds for $(e,f)$.
\medskip

{\it Proof of Claim 1.}
 Suppose $e=abc$, where $a,b,c$ are ordered as described earlier
  and suppose $f=acd$. Then (a) holds by $e$ being light and by \eqref{weight-bounds}.
(b) holds, since either $d(ab)\geq  \lfloor m/3\rfloor+2$ or $d(ab)= \lfloor m/3\rfloor+1$ and
$d'_{e,f}(a,b)=d(ab)-1$ (because $abd\notin G$). Similarly, (c) holds, since either $d(bc)\geq  \lfloor 2m/3\rfloor+2$ or 
$d(bc)= \lfloor 2m/3\rfloor+1$ and
$d'_{e,f}(b,c)=d(bc)-1$ (because $bcd\notin G$). 
Now, if $f$ is also a light edge then (d)  holds since $w(f)-\frac{3}{m}<0<\frac{3}{m}-w(e)$. \qed

By Claim 1, we may henceforth assume that every marked edge is heavy.
We will now use a discharging procedure to find our pair $(e,f)$.
Let the initial charge $ch(e)$ of every edge
$e$ in $G$ equal to $w(e)$. Then $\sum_{e\in G} ch(e)=\sum_{e\in G} w(e)=|\partial (G)|$.
We will redistribute charges among the edges of $G$ so that the total sum of charges
does not change and the resulting charge of each heavy edge remains at least $w_0$.  

The discharging rule is as follows.
 Suppose a heavy edge $f$ was marked by exactly $q$ light edges. If $q=0$, then let the new charge $ch^*(f)$
 equal $ch(f)$. Otherwise, let $f$ transfer to each light edge $e$ that marks it a charge  of $(ch(f)-w_0)/q$ so that $ch^*(f)=w_0$.
 It is easy to see that the total charge does not change in this discharging process.
 Hence, by \eqref{avg-wt}, there is an edge $e$ with $ch^*(e)<w_0$. By our discharging rule, $e$ must be a light edge.
Suppose $e$ marked $p$ edges. In each of Cases M1,M2, M3, $e$ marks at least one edge. So $p>0$.
Among all edges $e$ marked, let $f$ be one that gave the least charge to $e$. By definition, $f$ gave $e$ a charge of
at most $(ch^*(e)-ch(e))/p<(w_0-ch(e))/p$. We claim that the pair $(e,f)$ satisfies the lemma. Suppose $e=abc$, where $a,b,c$
are ordered as before, and suppose $f=acd$.
By Claim 1, (a), (b), and (c)  hold. It remains to prove (d). If all three pairs in $f$ are good, then $w(f)<\frac{3}{m}$,
contradicting $f$ being heavy. So, at most two of the pairs in $f$ are good. By our earlier discussion, $ac$ is good.
If one of $ad$ and $cd$ is also good, then the second part of (d) holds. So we may assume that $ac$ is the
only good pair in $f$.  Let $q$ be the number of the light edges that marked $f$. By the marking process,
a light edge only marks edges containing its high side and the high side is a good pair. Since $ac$ is the only
good pair in $f$, each of the $q$ light edges that marked $f$ contains $ac$ and has $ac$ as its high side.
 
First, suppose that  Case M1 was applied to $e$.
Then all the edges containing $ac$ other than $e$ were marked, which by our assumption must be heavy.
In particular, this implies that $q=1$. By our rule, $f$ gave $e$ a charge of $ch(f)-w_0$.  
By our choice of $f$, each  of the $d(ac)-1\geq m$ edges of $G$ containing $ac$ (other than $e$) gave $e$ a charge of at least $ch(f)-w_0$. 
Hence,
$w_0>ch^*(e)\geq ch(e)+m(ch(f)-w_0)$, from which the first part of (d) follows.

Next, suppose that Case M2 was applied to $e$. Then $d(ab)\leq \lfloor m/3 \rfloor +1$.
If $q>\lfloor m/3\rfloor+1$, then one of light edges containing $ac$, say $acx$, satisfies that
$abx\notin G$. By rule, $e$ marked $acx$, contradicting our assumption that no light edge was marked. So
$q\leq \lfloor m/3\rfloor+1$. Similarly if Case 3 was applied to $e$ then
$q\leq \lfloor 2m/3\rfloor+1$. In both of these cases,
$e$ marked at least $\frac{(m+3)(m-3)}{9}$ edges, and by the choice of $f$, each of these edges gave to $e$ charge at least $(ch(f)-w_0)/q$.
Since $ch^*(e)<w_0$, we conclude
$$w_0-ch(e)>\frac{(m+3)(m-3)}{9}\frac{ch(f)-w_0}{q}\geq \frac{(m+3)(m-3)}{3(2m+3)}(ch(f)-w_0).
$$
Since $m\geq 20$, this means
$$\frac{ch(f)-w_0}{w_0-ch(e)}<\frac{3(2m+3)}{(m+3)(m-3)}\leq \frac{3\cdot 45}{24\cdot 18}=\frac{5}{16}<\frac{1}{3}.
$$
So, the first part of (d) holds. \qed

For an edge $e$, by  $\mc(e)$ we denote the {\em minimum codegree} over all three pairs of vertices in $e$.

\begin{lem} \label{discharging2}
Let $G$ be a $3$-graph satisfying $e(G)> \gamma |\partial(G)|$. 
Let $w$ be the default weight function on $E(G)$ and $\partial(G)$.
Then there exists a pair of edges $e,f$ with $|e\cap f|=2$ such that 
\begin{enumerate}
\item $w(e)< \frac{1}{\gamma}$,
\item $d(e\cap f)=\mc(e)$,
\item $w(f)<\frac{1}{\gamma}+\frac{3}{\mc(e) -1}(\frac{1}{\gamma}-w(e))$.
\end{enumerate}
\end{lem}
\noindent
{\em Proof.}\enskip 
For convenience, let $w_0=\frac{1}{\gamma}$. As in the proof of Lemma \ref{discharging1},
call an edge $e$ with $w(e)<w_0$ {\it light} and an edge $e$ with $w(e)\geq w_0$ {\it heavy}. As before, the average average of $w(e)$ over all
$e$ is $|\partial(G)|/e(G)<w_0$. For each light edge $e$, let us mark a pair of vertices in that has codegree 
$\mc(e)$. If $e$ is a light edge with a marked pair $xy$ and $f$ is another light edge containing $xy$,
then our statements already hold. So we assume that no marked pair of any light edge 
lies in another light edge. Let us initially assign a charge of $w(e)$
to each edge $e$ in $G$. Then the average charge of an edge in $G$ is 
less than $w_0$.
 We now apply the following discharging rule.  For each heavy edge $f$,
transfer $\frac{1}{3}(w(f)-w_0)$ of the charge to each
light edge $e$ whose marked pair is contained in $f$.
Note that for each $f$ there are at most $3$ such $e$.
In particular, each heavy edge  still has charge at least $w_0$ after the discharging.

Since discharging does not change the total charge, there exists some edge
$e$ with charge less than $w_0$.  By the previous sentence, $e$ is a light edge in $G$ .
Let $xy$ be its marked pair. There are  $\mc(e)-1$ other edges 
containing it, each of which is heavy. Each such edge $f$ has given a charge of $\frac{1}{3}(w(f)-w_0)$ to $w_0$. For $e$ to still have a charge less than $w_0$,  one of these edges $f$ satisfies $\frac{1}{3}(w(f)-w_0)
< \frac{w_0-w(e)}{\mc(e)-1}$. Hence $w(f)<\frac{1}{\gamma}+\frac{3}{\mc(e)-1}(\frac{1}{\gamma}-w(e))$.\qed

Our third lemma proves a special case of Theorem ~\ref{two-central}.

\begin{lem} \label{two-central-weaker}
Let $T$ be a tight $3$-tree with $t\geq 5$ edges. Suppose $T$ has a trunk $\{e_1,e_2\}$ of size $2$ such that
$d_T(e_1\cap e_2)\geq \lfloor\frac{t-1}{3}\rfloor+2$. Let $G$ be an
$n$-vertex $3$-graph that does not contain $T$.
Then $e(G)\leq \frac{t-1}{3}|\partial(G)|$. 
\end{lem}


\noindent
{\em Proof.}\enskip 
  For convenience, let $m=t-1$.
Let $G$ be a $3$-graph with $e(G)>\frac{m}{3}|\partial(G)|$. Then $G$ contains
a subgraph $G'$ such that $e(G')> \frac{m}{3}|\partial (G')|$ and $\delta_2(G')> \frac{m}{3}$.
For convenience, we assume $G$ itself satisfies these two conditions. Let $w$ be the default weight
function on $E(G)$ and $\partial(G)$. Then $G$ satisfies the conditions of 
 Lemma \ref{discharging1}. Let the edges $e=abc$ and $f=adc$ satisfy the claim of that lemma, where
 $a,b,c$ are ordered as in Lemma \ref{discharging1}.
 In particular, by (a), $e$ is light and $ac$ is good, i.e. $d(ac)\geq m+1$. 
 By our assumptions, $d(ab)\leq d(bc)$. By parts (b) and (c),
   \begin{equation}\label{DL}
d'_{e,f}(a,b)\geq \left\lfloor\frac{m}{3}\right\rfloor\quad\mbox{\em and} \quad 
d'_{e,f}(c,b)\geq \left\lfloor\frac{2m}{3}\right\rfloor.
 \end{equation}  
 We rename pairs $\{a,d\}$ and $\{c,d\}$ as $D_1$ and $D_2$ so that $d'_{e,f}(D_1)=\min\{d'_{e,f}(a,d),d'_{e,f}(c,d)\}$ and
$d'_{e,f}(D_2)=\max\{d'_{e,f}(a,d),d'_{e,f}(c,d)\}$.
 We claim that in these terms,
   \begin{equation}\label{DR}
d'_1:=d'_{e,f}(D_1)\geq \left\lfloor\frac{m}{3}\right\rfloor-1\quad\mbox{\em and} \quad 
d'_2:=d'_{e,f}(D_2)\geq \left\lfloor\frac{m}{3}\right\rfloor.
 \end{equation}  
 By (\ref{d'}) and the fact that $\delta_2(G)>\frac{m}{3}$, $d'_1, d'_2\geq \lfloor \frac{m}{3}\rfloor -1$. 
 We will use part (d) of Lemma \ref{discharging1} to show that $d'_2\geq \lfloor \frac{m}{3}\rfloor$. If the second part of (d) holds, then
 $d'_2\geq m-1$ and we are done. So suppose the first part of Lemma 3.1 (d) holds instead, i.e. $3(w(f)-w_0)<(w_0-w(e))$. Then $w(f)<\frac{4}{3}w_0=\frac{4}{m}$. If
 $d'_1= d'_2=\left\lfloor\frac{m}{3}\right\rfloor-1$, then $d(D_1)=d(D_2)=\left\lfloor\frac{m}{3}\right\rfloor+1$ and hence
 $$w(f)>\frac{2}{\left\lfloor\frac{m}{3}\right\rfloor+1}\geq \frac{6}{m+3}\geq \frac{4}{m}
 $$
 when $m>9$, a contradiction. Thus, $d'_2\geq \lfloor \frac{m}{3} \rfloor$ and ~(\ref{DR}) holds.

By our assumption, $T$ has a trunk $\{e_1,e_2\}$ with $d_T(e_1\cap e_2)\geq \lfloor\frac{m}{3}\rfloor+2$. 
Suppose $e_1=xyu$ and $e_2=xyv$ so that
$e_1\cap e_2=xy$.  By our assumption, each edge in $E(T)\setminus \{e_1,e_2\}$
contains a pair in $e_1$ or $e_2$ and a vertex outside $e_1\cup e_2$.
For each  pair $B$ contained in $e_1$ or $e_2$, let
$N'_T(B)=N_T(B)\setminus \{x,y,u,v\}$ and $\mu(B)=|N'_T(B)|$.
Then $\mu(xy)=d_T(xy)-2$, and  $\mu(B)=d_T(B)-1$ for each $B\in \{xu,xv,yu,yv\}$,.
By definition, 
\begin{equation} \label{d-prime-sum}
\mu(xy)+\mu(xu)+\mu(xv)+\mu(yu)+\mu(yv)=t-2=m-1.
\end{equation}
Since $\mu(xy)=d_T(xy)-2\geq \lfloor\frac{m}{3}\rfloor>\frac{m}{3}-1$, we have
\begin{equation}  \label{sum}
\mu(xu)+\mu(xv)+\mu(yu)+\mu(yv)<\frac{2m}{3}.
\end{equation}

 We consider three cases, 
and in each case we find an embedding of $T$ into $G$.

\medskip

{\bf Case 1.} $d'_{e,f}(a,b)\geq \left\lfloor\frac{2m}{3}\right\rfloor$. Recall that by \eqref{DL},
$d'_{e,f} (c,b)\geq \lfloor \frac{2m}{3} \rfloor$.
By symmetry we may assume that $\mu(xu)+\mu(yu)\geq \mu(xv)+\mu(yv)$ and
that $\mu(xv)\geq \mu(yv)$. Then  by \eqref{sum}
$\mu(xv)+\mu(yv)\leq  \left \lfloor \frac{m}{3}\right \rfloor$, so we construct an embedding $\phi$ of $T$ into $G$ as follows.
 
 First, let $\phi(u)=b$  and $\phi(v)=d$. Then choose distinct $\phi(x),\phi(y)\in \{a,c\}$ so that $\phi(\{y,v\})=D_1$ and
 $\phi(\{x,v\})=D_2$.
This maps $e_1$ to $e$ and $e_2$ to $f$. Since $\mu(yv)<\frac{1}{4}\frac{2m}{3}=\frac{m}{6}$, by~(\ref{DR}) we can next
 map
$N'_T(yv)$ into $N'_G(D_1)$. Now, since  $\mu(yv)+\mu(xv)<\frac{1}{2}\frac{2m}{3}=\frac{m}{3}$, again by~(\ref{DR}) we can 
 map  $N'_T(xv)$  into $N'_G(D_2)\setminus \phi(N'_T(yv))$. If $\phi(x)=a, \phi(y)=c$, then
by the condition of Case 1 and~(\ref{sum}), we can map $N'_T(yu)$  into $N'_G(bc)\setminus \phi(N'_T(yv)\cup N'_T(xv)) $ and 
$N'_T(xu)$ into $N'_G(ac)\setminus \phi(N'_T(yv)\cup N'_T(xv))$. The case $\phi(x)=c, \phi(y)=a$ is similar.
Finally, embed $N'_T(xy)$ into
$N'_G(ac)$.

\medskip

{\bf Case 2.} $ \left\lfloor\frac{m}{3}\right\rfloor \leq d'_{e,f}(a,b)\leq \left\lfloor\frac{2m}{3}\right\rfloor-1$ and
$d'_1\geq \left\lfloor\frac{m}{3}\right\rfloor$.
Then we can strengthen the second part of~(\ref{DR}) to
 \begin{equation}\label{DR'}
d'_2\geq \left\lfloor\frac{m}{2}\right\rfloor.
 \end{equation} 
 Indeed, \eqref{DR} holds immediately if the second part of (d) holds in Lemma \ref{discharging1}; so we may assume
 $3(w(f)-w_0)<(w_0-w(e))$. By the condition of Case 2,
 $$w_0-w(e)\leq \frac{3}{m}-\frac{3}{2m+3}=\frac{3(m+3)}{m(2m+3)}.$$
 From this, we get 
 $$w(f)< \frac{3}{m}+\frac{(m+3)}{m(2m+3)}=\frac{7m+12}{m(2m+3)}.
 $$
 If $d'_2\leq \left\lfloor\frac{m}{2}\right\rfloor-1$, then
 $$w(f)>\frac{2}{d'_2+2}\geq 2\frac{2}{m+2},
 $$
 which is larger than $\frac{7m+12}{m(2m+3)}$ for $m\geq 24$. This contradiction proves~(\ref{DR'}).
  
 For convenience, suppose $D_1=cd$ (the case $D_1=ad$ is similar).  
By symmetry, we may assume that $\mu(xu)+\mu(yv)\leq \mu(yu)+\mu(xv)$ and 
that $\mu(yu)\geq \mu(xv)$. Then by \eqref{sum},
\begin{equation} \label{embed-upper2}
\mu(xu)+\mu(yv)\leq \left\lfloor\frac{m}{3}\right\rfloor, \quad \mu(xu)+\mu(yv)+\mu(xv)\leq \left\lfloor\frac{m}{2}\right\rfloor.
\end{equation}
We embed $T$ into $G$ by mapping  $x,y,u,v$ to $a,c,b,d$, respectively and
embedding in order $N'_T(yv)$ into $N'_G(cd)$, $N'_T(xu)$ into $N'_G(ab)$,
$N'_T(xv)$ into $N'_G(ad)$, $N'_T(yu)$ into $N'_G(bc)$, and $N'_T(xy)$ into
$N'_G(ac)$ greedily. 
Conditions \eqref{d-prime-sum},  \eqref{sum}, \eqref{DR'}  and \eqref{embed-upper2}
ensure that such an embedding exists.

\medskip
{\bf Case 3.} $ \left\lfloor\frac{m}{3}\right\rfloor \leq d'_{e,f}(a,b)\leq \left\lfloor\frac{2m}{3}\right\rfloor-1$ and
$d'_1= \left\lfloor\frac{m}{3}\right\rfloor-1$.
We now strengthen~(\ref{DR'}) to
 \begin{equation}\label{DR''}
d'_2\geq \left\lfloor\frac{2m}{3}\right\rfloor.
 \end{equation} 
 Indeed, exactly as in the proof of~(\ref{DR'}), we derive that $w(f)<\frac{7m+12}{m(2m+3)} $. 
If $d'_2\leq \lfloor  \frac{2m}{3}\rfloor-1$,  then
$$\frac{3}{2m+3}\leq\frac{1}{d'_2+2}<\frac{7m+12}{m(2m+3)}-\frac{1}{d'_1+2}\leq \frac{7m+12}{m(2m+3)}-\frac{3}{m+3},
$$
which is not true for $m\geq 20$. This proves~(\ref{DR''}).

As in Case 2,
suppose $D_1=cd$ (the case $D_1=ad$ is similar).  
By symmetry, we may assume that $\mu(xu)+\mu(yv)\leq \mu(yu)+\mu(xv)$ and 
that $\mu(xu)\geq \mu(yv)$. Then by \eqref{sum},
\begin{equation} \label{embed-upper3}
\mu(xu)+\mu(yv)\leq \left\lfloor\frac{m}{3}\right\rfloor, \quad \mu(yv)\leq \left\lfloor\frac{m}{6}\right\rfloor.
\end{equation}
We embed $T$ into $G$ by mapping  $x,y,u,v$ to $a,c,b,d$, respectively and
embedding in order $N'_T(yv)$ into $N'_G(cd)$, $N'_T(xu)$ into $N'_G(ab)$,
$N'_T(xv)$ into $N'_G(ad)$, $N'_T(yu)$ into $N'_G(bc)$, and $N'_T(xy)$ into
$N'_G(ac)$ greedily. 
Conditions  \eqref{d-prime-sum},  \eqref{sum}, \eqref{DR''} and \eqref{embed-upper3}
ensure that such an embedding exists.\qed


\section{Proof of Theorem \ref{two-central}}

We prove the shadow version of Theorem~\ref{two-central}, which immediately implies Theorem ~\ref{two-central}.

\medskip \noindent
{\bf Theorem~\ref{two-central}$'$.}
{\em
Let $t\geq 20$ be an integer. Let $T$ be a tight $3$-tree with $t$ edges and $c(T)\leq 2$.
If $G$ is an $r$-graph that
does not contain $T$ then $e(G)\leq \frac{t-1}{3} |\partial (G)|$.}

\begin{proof} First, let us point that in this proof, we exploit Lemma ~\ref{discharging2} and
will not need Lemma ~\ref{discharging1} in an explicit way.
Let $T$ be a tight $3$-tree with $t\geq 20$ edges that contains a trunk $\{e_1,e_2\}$ of size $2$.
For convenience, let $m=t-1$.
Let $G$ be a $3$-graph with $e(G)>\frac{m}{3}|\partial(G)|$. We
prove that $G$ contains $T$. As before we may assume that $\delta_2(G)>\frac{m}{3}$.
Let $w$ be the default weight function on $E(G)$ and $\partial(G)$.
By Lemma \ref{discharging2}, there exist  edges $e$ and $f$ in $G$ such that $d(e\cap f)=\mc(e)$,
$w(e)<\frac{3}{m}$, and
(using $m\geq 19$)
\begin{equation} \label{f-weight-upper}
\mbox{if $d(e\cap f)> \frac{m}{2}$, then} \; w(f)<\frac{3}{m}+\frac{3}{( m+1)/2-1}\left(\frac{3}{m}-w(e)\right)
\leq \frac{3}{m}+\frac{1}{3}\left (\frac{3}{m}-w(e)\right) \leq \frac{4}{m},
\end{equation}
and
\begin{equation} \label{f-weight-upper2}
\mbox{if $d(e\cap f)\leq \frac{m}{2}$, then} \; w(f)<\frac{3}{m}+\frac{3}{\lceil m/3\rceil-1}\left(\frac{3}{m}-w(e)\right)
\leq \frac{3}{m}+\frac{1}{2}\left (\frac{3}{m}-\frac{2}{m}\right) = \frac{4}{m}.
\end{equation}


Suppose $e=acb$ and $f=acd$, so that $e\cap f=ac$. 
For each pair $D$ contained in $e$ or $f$, let $N'_G(D)=N_G(D)\setminus \{a,b,c,d\}$
and $d'_G(D)=|N'_G(D)|$. Then $d'_G(D)\geq d_G(D)-2$.
Consider $T$. Suppose $e_1=xyu$ and $e_2=xyv$, so that $e_1\cap e_2 =xy$.
If $d_T(xy)\geq \lfloor\frac{m}{3}\rfloor+2$, then we  apply Lemma~\ref{two-central-weaker}
and are done. Hence we may assume that 
\[d_T(xy)\leq \left\lfloor \frac{m}{3}\right\rfloor+1.\]
For each pair $B$ contained in $e_1$ or $e_2$,
let $N'_T(B)=N_T(B)\setminus \{x,y,u,v\}$ and let $\mu(B)=|N'_T(B)|$. Then $\mu(xy)=d_T(xy)-2$
and $\mu(B)=d_T(B)-1$ for the other pairs. Also, we have
\begin{equation}\label{d-prime-sum2}
\mu(xu)+\mu(yu)+\mu(xv)+\mu(yv)+\mu(xy)=m-1.
\end{equation}
Since $\mu(xy)=d_T(xy)-2\leq \frac{m}{3}-1$, 
\begin{equation} \label{d-prime-sum3}
  \mu(xy)+\frac{i}{4}(m-1-\mu(xy))\leq \frac{m}{3}+\frac{im}{6}-1\qquad\qquad \forall i\in [4].
\end{equation}

Let us view $e,f$ as glued together at $ac$ with $e$ on the left and
$f$ on the right. Let
\begin{eqnarray*}
L_{max}=\max\{d_G(ab), d_G(bc)\}, && \quad L_{min}=\min\{d_G(ab), d_G(bc)\},\\
R_{max}=\max\{d_G(ad), d_G(cd)\},&& \quad R_{min}=\min\{d_G(ad), d_G(cd)\}.
\end{eqnarray*}

Since $d(ac)=\mc(e)$, 
 $L_{max}\geq L_{min} \geq d_G(ac)$. Since $w(e)<\frac{3}{m}$, we have 

\begin{equation} \label{Lmax}
L_{max}>m.
\end{equation}

We consider two cases. In each case, we find an embedding of $T$ into $G$.

\medskip

{\bf Case 1.} $L_{min}>m$. 
This implies
$d'_G(ab), d'_G(bc)\geq m-1$.
By symmetry, we may assume that $d_G(ad)\geq d_G(cd)$
so that  $d_G(ad)=R_{max}$ and $d_G(cd)=R_{min}$.
Now, consider $T$. By symmetry, we may assume
that $\mu(xu)+\mu(yu)\geq \mu(xv)+\mu(yv)$ and that $\mu(xv)\geq \mu(yv)$.
Then $\mu(yv)\leq \frac{1}{4}(m-1-\mu(xy))$ and $\mu(xv)+\mu(yv)\leq \frac{1}{2}{}(m-1-\mu(xy))$.
This, together with \eqref{d-prime-sum3} implies
\begin{eqnarray}\label{new-embed-upper1a}
&\mu(yv)\leq \left\lfloor \frac{m}{4}\right \rfloor, \quad &\mu(xv)+\mu(yv)\leq \left\lfloor\frac{m}{2}\right\rfloor-1, \nonumber\\
&\mu(yv)+\mu(xy)\leq \left\lfloor\frac{m}{2}\right\rfloor-1, \quad
&\mu(xv)+\mu(yv)+\mu(xy)\leq \left\lfloor\frac{2m}{3}\right\rfloor-1.
\end{eqnarray}


\medskip

{\bf Case 1.1.} $d_G(ac)>\frac{2m}{3}$.
By \eqref{f-weight-upper}, $\frac{1}{R_{max}}+\frac{1}{R_{min}}<w(f)<\frac{4}{m}$, so
$R_{max}>\frac{m}{2}$. Since $\delta_2(G)>\frac{m}{3}$, we have $R_{min}>\frac{m}{3}$. 
Hence
\begin{equation} \label{new-embed-lower1a}
d'_G(ab), d'_G(bc)\geq m-1, \quad d'_G(ac)\geq \left\lfloor \frac{2m}{3}\right\rfloor-1, \quad  d'_G(ad)\geq \left\lfloor \frac{m}{2}\right\rfloor-1,\quad 
d'_G(cd)\geq \left\lfloor \frac{m}{3}\right \rfloor -1.
\end{equation}
Now we can embed $T$ into $G$ as follows. 
First, we map $x,y,u,v$ to $a,b,c,d$ respectively. This maps $e_1$ to $e$ and $e_2$ to $f$. 
Then we map $N'_T(yv)$ into $N'_G(cd)$ followed by $N'_T(xv)$ into $N'_G(ad)$. 
Next, we map $N'_T(xy)$ into $N'_G(ac)$, $N'_T(yu)$ into $N'_G(bc)$, and $N'_T(xu)$ into $N'_G(ab)$
in that order. Conditions \eqref{new-embed-upper1a} and \eqref{new-embed-lower1a} ensure that such an embedding
exists.

\medskip

{\bf Case 1.2.} $d_G(ac)\leq \frac{2m}{3}$. Then $w(e)\geq \frac{3}{2m}$.
If  $d_G(ac)> \frac{m}{2}$, then by \eqref{f-weight-upper}, $w(f)<\frac{3}{m}+\frac{1}{3}(\frac{3}{m}-\frac{3}{2m})=\frac{7}{2m}$.
On the other hand, if $d_G(ac)\leq \frac{m}{2}$, then $w(e)\geq \frac{2}{m}$ and 
by \eqref{f-weight-upper2}, $w(f)<\frac{3}{m}+\frac{1}{2}(\frac{3}{m}-\frac{2}{m})=\frac{7}{2m}$.
So in any case,
$$\frac{1}{R_{max}}+\frac{1}{R_{min}}<w(f)-w(ac)<\frac{7}{2m}-\frac{3}{2m}=\frac{2}{m}.$$

Then $R_{max}>m$ and $R_{min}>\frac{m}{2}$. Also, since $\delta_2(G)>\frac{m}{3}$, we have $d_G(ac)>\frac{m}{3}$. Hence,
\begin{equation}\label{new-embed-lower1b}
d'_G(ab), d'_G(bc)\geq m-1,  \quad d'_G(ac)\geq \left\lfloor \frac{m}{3}\right\rfloor -1, \quad d'_G(ad)\geq m-1,\quad d'_G(cd)\geq \left\lfloor \frac{m}{2}\right \rfloor-1.
\end{equation}
Now we can embed $T$ into $G$ as follows. 
First, we map $x,y,u,v$ to $a,b,c,d$ respectively. This maps $e_1$ to $e$ and $e_2$ to $f$. 
Then we map $N'_T(xy)$ into $N'_G(ac)$. This is doable since $d'_T(xy)=d_T(xy)-2\leq \lfloor \frac{m}{3}\rfloor-1$
while $d'_G(ac)\geq \lfloor \frac{m}{3}\rfloor-1$.
Then we map $N'_T(yv)$ into $N'_G(cd)$ followed by $N'_T(xv)$ into $N'_G(ad)$. 
Next, we map $N'_T(yu)$ into $N'_G(bc)$, and $N'_T(xu)$ into $N'_G(ab)$
in that order. Conditions \eqref{new-embed-upper1a} and \eqref{new-embed-lower1b} ensure that such an embedding
exists.

\medskip

{\bf Case 2.} $L_{min}\leq m$.
By symmetry, we may assume that $d_G(ab)\geq d_G(bc)$ so that $d_G(ab)=L_{max}$ and $d_G(bc)=L_{min}$.
We have $\frac{1}{L_{min}}+\frac{1}{d_G(ac)}<w(e)<\frac{3}{m}$.
Since $d(ac)=\mc(e)$, 
$d_G(ac)\leq L_{min}\leq m$. This yields  $L_{min}>\frac{2m}{3}, \frac{m}{2}<d_G(ac)\leq m$, and 
$w(e)>\frac{2}{m}$.  By \eqref{Lmax}, $L_{max}>m$.  Thus, 
\begin{equation} \label{new-embed-lower2a}
d'_G(ab)\geq m-1, \quad d'_G(bc)\geq \left\lfloor \frac{2m}{3}\right\rfloor-1, \quad d'_G(ac)\geq \left \lfloor \frac{m}{2}\right\rfloor-1.
\end{equation}

Since $d_G(ac)>m/2$, by \eqref{f-weight-upper}, 
\begin{equation} \label{Rmaxmin}
w(f)<\frac{3}{m}+\frac{1}{3}\frac{1}{m}=\frac{10}{3m} \mbox{ and }
\frac{1}{R_{max}}+\frac{1}{R_{min}}\leq w(f)-\frac{1}{d_G(ac)}<\frac{10}{3m}-\frac{1}{m}=\frac{7}{3m}.
\end{equation}

{\bf Case 2.1} $R_{max}>m$.
%
%
 By our assumption and \eqref{Rmaxmin}, 
\[R_{max}>m, \quad R_{min}>\frac{3m}{7}.\]

\medskip

First suppose that $d_G(ad)\geq d_G(cd)$. Then
\begin{equation}\label{new-embed-lower2b}
d'_G(ad)\geq m -1,\quad d'_G(cd)\geq \left \lfloor \frac{3m}{7}\right\rfloor-1.
\end{equation}

\medskip

By symmetry, we may assume that 
$\mu(xu)+\mu(xv)\geq \mu(yu)+\mu(yv)$ and that $\mu(yu)\geq \mu(yv)$.
Then by these assumptions and \eqref{d-prime-sum3}, we have
\begin{eqnarray}\label{new-embed-upper2a}
\mu(yv)\leq \left\lfloor \frac{m}{4}\right\rfloor-1,\quad \mu(yv)+\mu(xy)\leq \left\lfloor \frac{m}{2}\right\rfloor-1,
\quad \mu(yv)+\mu(xy)+\mu(yu)\leq \left\lfloor \frac{2m}{3}\right\rfloor -1.
\end{eqnarray}

Now we can embed $T$ into $G$ as follows. 
First, we map $x,y,u,v$ to $a,b,c,d$ respectively. This maps $e_1$ to $e$ and $e_2$ to $f$. 
Then we map $N'_T(yv)$ into $N'_G(cd)$ followed by $N'_T(xy)$ into $N'_G(ac)$. 
Next, we map $N'_T(yu)$ into $N'_G(bc)$, $N'_T(xv)$ into $N'_G(ad)$, and $N'_T(xu)$ into $N'_G(ab)$
in that order. Conditions \eqref{new-embed-lower2a}, \eqref{new-embed-lower2b} and \eqref{new-embed-upper2a} 
ensure that such an embedding exists.

Next, suppose that $d_G(cd)\geq d_G(ad)$. Then
\begin{equation}\label{new-embed-lower2c}
d'_G(ad)\geq  \left \lfloor \frac{3m}{7}\right\rfloor-1,\quad d'_G(cd)\geq m-1.
\end{equation}

By symmetry, we may assume that 
$\mu(xu)+\mu(yv)\geq \mu(xv)+\mu(yu)$ and that $\mu(yu)\geq \mu(xv)$.
By these assumptions and \eqref{d-prime-sum3}, we have
\begin{eqnarray}\label{new-embed-upper2c}
\mu(xv)\leq \left\lfloor \frac{m}{4}\right\rfloor-1,\quad \mu(xv)+\mu(xy)\leq \left\lfloor \frac{m}{2}\right\rfloor-1,
\quad \mu(xv)+\mu(xy)+\mu(yu)\leq \left\lfloor \frac{2m}{3}\right\rfloor -1.
\end{eqnarray}

Now we can embed $T$ into $G$ as follows. 
First, we map $x,y,u,v$ to $a,b,c,d$ respectively. This maps $e_1$ to $e$ and $e_2$ to $f$. 
Then we map $N'_T(xv)$ into $N'_G(ad)$ followed by $N'_T(xy)$ into $N'_G(ac)$. 
Next, we map $N'_T(yu)$ into $N'_G(bc)$, $N'_T(yv)$ into $N'_G(cd)$, and $N'_T(xu)$ into $N'_G(ab)$
in that order. Conditions \eqref{new-embed-lower2a}, \eqref{new-embed-lower2c} and \eqref{new-embed-upper2c} 
ensure that such an embedding exists.

\medskip

{\bf Case 2.2} $R_{max}\leq m$.
%
%
Since $R_{min}\leq R_{max}\leq m$, by \eqref{Rmaxmin}, we again have
 $R_{max}>\frac{6m}{7}$, and 
 $$\frac{1}{R_{min}}< \frac{7}{3m}-\frac{1}{m}=\frac{4}{3m};\quad\mbox{ so     }\quad R_{min}   >   \frac{3m}{4}.$$
By \eqref{Rmaxmin}, $w(f)<\frac{10}{3m}$. Also, $\frac{1}{L_{min}}\geq \frac{1}{L_{max}}\geq \frac{1}{m}$. Hence,
 \begin{equation} \label{ac-new}
w(ac)<\frac{10}{3m}-\frac{2}{m}=\frac{4}{3m} \quad\mbox{\em and hence}\quad  d'_G(ac) \geq \left\lfloor \frac{3m}{4}\right\rfloor-1.
\end{equation}

First, suppose that $d_G(ad)\geq d_G(cd)$. Then
\begin{equation} \label{new-embed-lower-2d}
d'(ad)\geq \left\lfloor \frac{6m}{7}  \right \rfloor -1, \quad d'(cd)\geq \left\lfloor \frac{3m}{4} \right\rfloor -1. 
\end{equation} 
By symmetry, we may assume that $\mu(xu)+\mu(xv)\geq \mu(yu)+\mu(yv)$ and that $\mu(xu)\geq \mu(xv)$. In particular,
 \begin{equation}\label{xu}
\mu(xu)\geq \frac{1}{4}\left(m-1-\mu(xy)\right)\geq \frac{1}{4}\left(m-1-\frac{m}{3}+1\right)=\frac{m}{6}.
 \end{equation}  
By \eqref{d-prime-sum3}, \eqref{new-embed-lower2a}, \eqref{new-embed-lower-2d}, and  (\ref{xu}), we can greedily embed $T$ into $G$
by mapping $x,y,u,v$ to $a,c,b,d$, respectively and mapping in order $N'_T(yv)$ into $N'_G(cd)$, $N'_T(xy)$ into $N'_G(ac)$,
$N'_T(yu)$ into $N'_G(bc)$, $N'_T(xv)$ into $N'_G(ad)$, and $N'_T(xu)$ into $N'_G(ab)$.

Next, suppose that $d_G(cd)\geq d_G(ad)$. Then $d'(ad)\geq \left\lfloor \frac{3m}{4}  \right \rfloor -1$ and  $d'(cd)\geq \left\lfloor \frac{6m}{7} \right\rfloor -1$. 
By symmetry, we may assume that $\mu(xu)+\mu(yv)\geq \mu(xv)+\mu(yu)$ and that $\mu(xu)\geq \mu(yv)$.
Again,~(\ref{xu}) holds.
We can greedily embed $T$ into $G$
by mapping $x,y,u,v$ to $a,c,b,d$, respectively and mapping in order $N'_T(yu)$ into $N'_G(bc)$, $N'_T(xy)$ into $N'_G(ac)$,
$N'_T(xv)$ into $N'_G(ad)$, $N'_T(yv)$ into $N'_G(cd)$, and $N'_T(xu)$ into $N'_G(ab)$.
\end{proof}

\paragraph{Acknowledgement.}
This research was partly conducted during an American Institute of Mathematics Structured Quartet Research Ensembles workshop,
and we gratefully acknowledge the support of AIM.

{\small

\begin{tabular}{ll}
\begin{tabular}{l}
{\sc Zolt\'an F\" uredi} \\
Alfr\' ed R\' enyi Institute of Mathematics \\
Hungarian Academy of Sciences \\
Re\'{a}ltanoda utca 13-15 \\
H-1053, Budapest, Hungary \\
E-mail:  \texttt{zfuredi@gmail.com}.
\end{tabular}
& \begin{tabular}{l}
{\sc Tao Jiang} \\
Department of Mathematics \\ Miami University \\ Oxford, OH 45056, USA. \\ E-mail: \texttt{jiangt@miamioh.edu}. \\
$\mbox{ }$
\end{tabular} \\ \\
\begin{tabular}{l}
{\sc Alexandr Kostochka} \\
University of Illinois at Urbana--Champaign \\
Urbana, IL 61801 \\
and Sobolev Institute of Mathematics \\
Novosibirsk 630090, Russia. \\
E-mail: \texttt {kostochk@math.uiuc.edu}.
\end{tabular} & \begin{tabular}{l}
{\sc Dhruv Mubayi} \\
Department of Mathematics, Statistics \\
and Computer Science \\
University of Illinois at Chicago \\
Chicago, IL 60607. \\
\texttt{E-mail: mubayi@uic.edu}.
\end{tabular} \\ \\
\begin{tabular}{l}
{\sc Jacques Verstra\"ete} \\
Department of Mathematics \\
University of California at San Diego \\
9500 Gilman Drive, La Jolla, California 92093-0112, USA. \\
E-mail: {\tt jverstra@math.ucsd.edu.}
\end{tabular}
\end{tabular}
}

\end{document}